**Carlos Grandjot**
**Drei Jahrzehnte der Mathematik in Chile: 1930 – 1960**
Von Claudio Gutiérrez und Flavio Gutiérrez
Übersetzt von Karin Erxmeyer

## 1 Einleitung

Sieht man sich die Entwicklung der Wissenschaften in Chile genauer an, besonders die der Mathematik, so kommt man zu dem Ergebnis, dass ein langer und steiniger Weg abgeschritten werden musste, um an ihrer jetzigen Reife anzukommen. Viele Personen trugen zu dieser Entwicklung bei, davon auch einige Ausländer, die Chile zu ihrer zweiten Heimat machten. Carlos Grandjot war einer von ihnen. Um seinen Beitrag zur Mathematik zu verstehen, ist es notwendig, sich einen Überblick über die Entwicklung dieser Disziplin zu verschaffen.

Ihre Anfänge gehen zurück auf die letzten Jahre der chilenischen Kolonialgeschichte, als einerseits 1758 der Mathematik-Lehrstuhl an die Königliche Universität von San Felipe vergeben wurde und andererseits 1799 mathematische Vorlesungen an der Akademie von San Luis begannen. Bei Beginn der Unabhängigkeit vereinigten sich diese beiden Institutionen zum Nationalen Institut im Jahre 1813, aus dem 1824 die ersten Landvermesser der Republik mit abgeschlossener Ausbildung hervorgingen. Die Inhalte dieser Kurse wurden von den Professoren selbst festgelegt und anhand eigener Aufzeichnungen vorgetragen. Der spanische Ingenieur Andrés Antonio Gorbea, der 1826 nach Chile kam, organisierte die mathematische Lehre dann nach europäischem Stil neu. Er entwarf ein Programm, das 1831 offiziell anerkannt wurde, und legte seine Ideen für seine Vorlesungen in der polytechnischen Schule schriftlich nieder. Die Leitlinie dieser Institution war, die Mathematik als anwendungsbezogene Wissenschaft (*ciencia „útil"*) zu etablieren. In diesem Geiste verblieben die mathematischen Wissenschaften in Chile bis zum Ende des 19. Jahrhunderts. Der Richtungsgeber für deren Entwicklung war die Ingenieurschule.

Mit der Ankunft der deutschen Professoren in 1889, die das Pädagogische Institut aufbauen sollten, erreichten die mathematischen Wissenschaften eine neue Dimension. Tafelmacher und Poenisch, beide promoviert in den exakten Wissenschaften, bauten die Mathematik nun um in eine autonome und kulturelle Disziplin, die als Kernstück allen Wissens der Lehre dienen sollte, eine Eigenschaft, die sie in Europa schon seit Beginn des Jahrhunderts hatte. Sie entwarfen Lehrprogramme und Unterrichtstexte und veröffentlichten Artikel über aktuelle Themen in den Annalen der Universität von Chile. Dieser neuartige Auftritt der Mathematik brachte unter einigen Intellektuellen weitreichende Veränderungen hervor. Valentín Letelier kommentierte 1895 in seinem Buch „Der Kampf um die Kultur": „Seit so vielen Jahren ist in Chile nicht mehr über mathematische Themen geschrieben worden, dass die letzten Schülergenerationen im Glauben erzogen wurden, dass diese Wissenschaft vollständig mumifiziert und nicht zu einer größeren Entwicklung fähig sei."

Die Impulse von Poenisch und Tafelmacher konzentrierten sich besonders auf die Lehre. Dabei bildeten sie eine ganze Legion von treuen Schülern aus, die sich über ganz Chile ausbreiteten. Der bereits pensionierte Poenisch konnte daher 1929 folgendermaßen seiner Zufriedenheit Ausdruck verleihen: „Die Lehre dieses Faches befindet sich vertrauensvoll in den Händen von Personen, die ihre Aufgaben kennen und erfüllen können (…). Ihre Bereitschaft und Arbeitsauffassung lassen mich in Ruhe ausruhen." Die Zufriedenheit des Meisters war sicherlich berechtigt.

Aber nur zwei notwendige Aspekte einer reifen Wissenschaft wurden damals in der Mathematik in Chile bereits ausreichend berücksichtigt: Sie war sowohl anwendungsbezogen als auch eine Wissens-Sammlung. Man musste jedoch feststellen, dass noch ein dritter Aspekt fehlte, um ihre Weiterentwicklung zu stimulieren und ihr Stocken zu verhindern. Dieser fehlende Aspekt ist die *Wissenschaft als kreativer Prozess.* Carlos Grandjot übernahm die Eingliederung dieses dritten Aspektes in die exakten Wissenschaften – als Pionier in Chile.

Seit seiner Ankunft im Land beteiligte sich Grandjot an der Gründung von Institutionen für den Fortschritt und die Kultur der nationalen Wissenschaften. Er war Gründungsprofessor des Chilenischen Institutes (1930), der ersten Mathematischen Vereinigung Chiles (1953), dessen erster Präsident er wurde, des Institutes der mathematischen Forschung (1957) und weiterer Einrichtungen. Von seinem Lehrstuhl aus verbreitete und förderte er die fortschrittlichste Mathematik und moderne Physik unter seinen Schülern, manchmal mit Erfolg, manchmal auch ohne diesen, aber immer mit viel Enthusiasmus. Er arbeitete eng mit den akademischen Speerspitzen der wissenschaftlichen Forschung in den Institutionen zusammen und hatte die erforderlichen Kontakte mit den her-ausragenden internationalen Zentren, um die Forschung auch in Chile voranzubringen.

Sollten die Kenntnis seines Lebensweges und seines Lebenswerkes die Entwicklung der Mathematik in Chile in der damaligen Zeit erhellen können, dann hätte dieser Artikel seinen Zweck erfüllt.

Dr. Grandjot, wie ihn seine Schüler nannten, kam am 10. Mai 1929 aus Deutschland in Chile an. Von der Regierung war er unter Vertrag genommen worden, um als Mathematikprofessor in den Lehranstalten der Republik seine Dienste zu absolvieren. Es waren wöchentlich fünfzehn Stunden Vorlesung mit den dazugehörigen Seminaren festgelegt.[1] Gleich nach seiner Ankunft begann er mit Vorlesungen der elementaren und höheren Mathematik, der Philosophie und der Physik im Pädagogischen Institut von Chile.[2] Der Vertrag hatte eine Laufzeit von zwei Jahren ab dem 9. April 1929.[3] Immer wieder verlängerte er diesen Vertrag, bis ihn in Chile der Beginn des Zweiten Weltkrieges überraschte. Dies trug mit zur Besiegelung seines Schicksals bei, sich zu entscheiden, in Chile zu bleiben und die Staatsbürgerschaft anzunehmen, was zusätzlich zum Pädagogischen Institut drei weiteren universitären Schulen erlaubte, sich seiner Dienste zu erfreuen: die Ingenieursschule der katholischen erzbischöflichen Universität Chile seit 1933, die Ingenieursschule der Universität Chile mit ihrem Kurs „Ergänzungen der höheren Mathematik" ab 1945 und die Architekturschule der oben erwähnten katholischen Universität ab dem Jahr 1953.

Die Aktivitäten von Karl Grandjot beschränkten sich nicht nur auf die Vorlesungen. Seine umfangreiche Bildung, seine besondere Gabe für die wissenschaftliche Forschung und sein Enthusiasmus ließen ihn auch an anderen Vereinigungen verschiedener Richtungen teilnehmen: er war Präsident der Chilenischen Vereinigung für Naturgeschichte, in der er an verschiedenen Originalarbeiten teilnahm; Mitglied des Rates des Deutsch-Chilenischen Bundes; Gesellschafter und Direktor der Vereinigung der Musik Mozarts in Santiago; in seiner Eigenschaft als Gesellschafter des Deutschen Wanderclubs durchwanderte er Chile in all seiner Ausdehnung; er bereiste Bolivien und Hochperu; besuchte mehrmals die Region Araucanien und sein ganzer persönlicher Stolz war, dass er das Mapudungún (die Sprache der Mapuche) direkt von den Indios des Südens gelernt hatte. Seine sprachlichen Fähigkeiten waren überragend: abgesehen von Deutsch sprach er Spanisch, Englisch, Französisch, Russisch, Mapudungún und Portugiesisch, zusätzlich zu Latein, Griechisch und Holländisch. Auf seinen Reisen interessierte er sich außerdem für Aymará und Quetschua. Er kannte sich in mehreren deutschen Dialekten aus und hatte Grundkenntnisse in Japanisch und Chinesisch, deren Schriftzeichen er mit Freude verglich.

Aus Deutschland kommend, eilte Grandjot ein verdienter Ruf als Wunderstudent voraus. Dies war vor allem darauf zurückzuführen, dass er mit Landau ein Axiomensystem der natürlichen Zahlen verbessert hatte. Diese Kritik an den Peano-Axiomen war als „Einwand von Grandjot" bekannt geworden. Außerdem konnte er viele Publikationen zu den schwierigsten mathematischen Gebieten, wie der Zahlentheorie, in den angesehensten Zeitschriften Europas vorweisen. Unter den Studenten des Pädagogischen Institutes der 40er Jahre kursierte das nicht bestätigte, aber auch nicht dementierte Gerücht, dass die chilenische Regierung Grandjot in einem strengen Wettbewerb unter 500 Anwärtern auf den Lehrstuhl ausgewählt hatte.[4]

Dr. Grandjot hatte aus jedwedem Blickwinkel eine erstaunliche Leichtigkeit um seine Vorlesungen einzuleiten. Er begann mit dem Aufstellen einiger Anfangsbehauptungen und im Folgenden zog er daraus die jeweils angemessensten Schlussfolgerungen für seine Zwecke.[5] Er verglich die Mathematik mit dem Bau eines Hauses. Als erstes kommen die Fundamente und erst danach die Errichtung des Gebäudes. Die Fundamente können entsprechend dem Untergrund breit sein und mit Überfluss an Material erstellt oder sehr schmal mit geringstem materiellen Aufwand. In Analogie kann man bei der Konstruktion einer mathematischen Theorie ein Axiomensystem mit einer großen, mehr als nötigen Zahl an Axiomen postulieren oder aber eines mit der genau richtigen Anzahl. Beide Optionen sind legitim, je nachdem, welche Strenge man an deren Ausarbeitung legt. Zu diesem Thema sein pädagogischer Rat: Die Intuition ist im schulischen Unterricht ein wertvolles Instrument, wenn man sie gut und mit Klugheit anwendet. Übermäßige Strenge oder Formalismus führen nämlich dazu, dass die Schüler die Lust an der Mathematik verlie-ren. Man muss dabei mit großer Vorsicht vorgehen.

---

[1] Decreto Supremo Nr. 1764, als Nachweis des Vertrages vom 17. Mai 1929 zwischen der Regierung und Herrn Karl Grandjot. Bericht des Universitären Rates, Regierungserlasse
[2] Lebenslauf [20]
[3] Seine Tochter Sigrid erinnerte sich jedoch, immer von ihrem Vater gehört zu haben, dass der Vertrag auf 5 Jahre abgeschlossen war, verlängerbar um jeweils weitere 5 Jahre.

[4] „Wir als seine Schüler zweifelten nicht im Geringsten am Wahrheitsgehalt dieser Anekdote." F. Gutiérrez, persönliche Erinnerungen, 1951.
[5] Dies wird bestätigt durch die Mitteilungen von Dr. Benedicto Chuaqui nach einer gemeinsamen Durchsicht aller privaten Aufzeichnungen über seine Klassen, die Grandjot ihm und Rolando Chuaqui überlassen hatte.

Zu seinem großen Können in der Forschung gesellte sich bei Grandjot das richtige Gespür für das Unterrichten. Stellte er bei seinen Schülern Müdigkeit oder geistige Erschöpfung fest, so warf er einen Witz oder eine passende Anekdote ein. So erzählte er nebenbei zur Theorie der Farben, dass er kürzlich eine rote Krawatte getragen hätte im Lehrerzimmer der Ingenieursschule Abraham Pérez, Professor der Geometrie und der Algebra des Pädagogischen Institutes, getroffen hätte. Er hätte sich ihm genähert und seine Krawatte mit den Worten gezeigt: „Wie gefällt dir meine Krawatte, Abraham?" Pérez, der sofort eine Hinterlist seines Kollegen vermutete, befühlte sie mit seinen Händen und antwortete lakonisch: „Es erscheint mir eine gute Qualität!", eine Antwort, die Grandjot mit Lachen aus vollem Halse wiedergab. Erst danach fügte er hinzu, dass Pérez an Farbenblindheit litt, so dass sich dann die Klasse ebenso über den Witz amüsieren konnte.[6] Hatte Grandjot so die Aufmerksamkeit der Klasse wieder hergestellt, nahm er seinen Unterricht wieder auf.

Zu seinem großen Humor[7] kam zudem die Schnelligkeit des Denkens. Diese zeigte er bei einer anderen Gelegenheit, nach dem Ende seiner Konferenz „Was ist das Leben?"[8]. Einer der Teilnehmer zog seine Taschenuhr aus der Tasche und sagte zu den gerade vorgetragenen Theorien hinterlistig: „Ich glaube, Herr Professor, meine Uhr besitzt Leben! Was meinen Sie dazu?" Schlagfertig antwortete Grandjot: „Sobald Sie, mein Herr, eine Abhandlung über die Fortpflanzung der Uhren vorgelegt haben, werde ich meine Meinung darüber mitteilen!"

Das Leben Grandjots war arbeitsreich und es fehlten auch nicht schwierige Zeiten. Aber er ging es trotzdem entspannt an. Seine größte Leidenschaft war das Bergsteigen auf dem Rücken eines Mulis und mit Zelt. Fast immer wurde er dabei von seiner Ehefrau Gertrudis Fritsche begleitet, die er in Göttingen 1926 geheiratet hatte.[9] Beide waren Liebhaber der Na-tur und botanische Experten. Auf ihren Exkursionen sammelten sie einheimische Pflanzen. Wie wir ihren Publikationen in der Zeitschrift der Naturgeschichte und in der Zeitschrift der Wissenschaftlichen Vereinigung Deutschlands entnehmen können, haben sie viele dieser Arten klassifiziert. Ein weiterer Zeitvertrieb von Grandjot war die Musik. Für besondere Gelegenheiten, wie Geburtstage, Taufen oder Hochzeiten komponierte er selbst Musik für seine Freunde. Die Partitur dazu erarbeitete er an seinem Schreibtisch und erst später überprüfte er sie am Klavier.[10]

Er liebte es, Blockflöte zu spielen und in seinem Haus organisierte er Treffen, an denen deutsche Studentengruppen zusammen mit ihm musizierten. Er sang regelmäßig im Chor des Dirigenten Jan Sparwaater und spielte in Aufführungen der Theatergruppe des Deutschen Reinhold Olszewski und im Laientheater der Familie von Viesling in Las Condes, wenn es sich ergab auch im Marionettentheater der gleichen Familie. Grandjot hatte viele deutsche und chilenische Freunde und mit einigen ging er an den Sonntagen auf Exkursionen. Als guter Kartenspieler hatte er mit zwei Namensvettern viel Spaß beim Skatspielen.[11]

Karl Grandjot Reins wurde am 23. August 1900 in Frankenberg, Deutschland, geboren. Er war der älteste Sohn des Postinspektors Konrad Grandjot Blume und von Luise Reins Remhof. Er hatte 2 Brüder: Erich, Straßenbauingenieur und Walter, der jüngste, Doktor der Physik, auf Akustik spezialisiert. Von seinen Schülerjahren erzählt Grandjot in seiner kurzen Autobiografie, dass er die Volksschule und die Oberrealschule in Kassel besuchte. Später, im Alter von 19 Jahren, schrieb er sich an der Universität von Göttingen ein, wo er rein und angewandte Mathematik, experimentelle und theoretische Physik und Philosophie studierte. Unter seinen wichtigsten Professoren waren so Angesehene wie die Mathematiker Edmund Landau, Richard Courant, David Hilbert und die Physiker Peter Debye und Max Born, beide mit dem Nobelpreis ausgezeichnet. Am 14. Februar 1922 promovierte er und wurde zur gleichen Zeit wissenschaftlicher Assistent Edmund Landaus. 1925 graduierte er zum Privatdozenten, ein Titel, der ihm die universitäre Lehre ermöglichte. Es folgte eine fruchtbare Zusammenarbeit mit Landau und eine Vielzahl von Publikationen in verschiedenen europäischen mathematischen Zeitschriften in den Jahren 1922 bis 1929. Er nahm in diesen Jahren an mehreren wissen-

---

[6] Erinnerungen des F. Gutiérrez, Schüler von Grandjot im Jahre 1951 im Kurs der Theoretischen Physik.
[7] Eine weitere Anekdote, die Autoren von B. Chuaqui erzählt: Immer wenn er in seinem Haus einen falsch verbundenen Telefonanruf erhielt, antwortete er spitzbübisch: „Sie haben nur drei Richtige!"
[8] Im Jahr 1947. Publiziert in „Impulso", der Zeitschrift des Zentrums der Ingeniere der Katholischen Universität. Siehe Literaturverzeichnis.
[9] Gertrudis kam im September 1929 fünf Monate nach Karl in Chile an, zusammen mit der gemeinsamen Tochter Sigrid, die im Februar 1929 in Paris geboren wurde, wenige Monate bevor Karl nach Chile ging. Sigrid studierte an der Universität von

Chile und erhielt 1953 den Abschluss als Mathematik- und Physiklehrerin.
[10] Erinnerungen seiner Tochter Sigrid. Brief an die Autoren.
[11] Skat ist das deutsche Nationalkartenspiel für 3 Mitspieler, entstanden 1810.

schaftlichen Kongressen teil, unter anderem am Internationalen Mathematikerkongress in Bologna im August 1928. Von 1928 bis 1929 erhielt er ein Stipendium der Rockefellerstiftung, um seine Studien in Paris zu intensivieren. Mitten in diesen Studien erreichte ihn der Ruf der chilenischen Regierung. Im April 1929 schiffte er sich nach Santiago ein und erreichte die chilenische Küste am 10. Mai desselben Jahres, wie wir am Anfang bereits erläutert haben.

## 2 Die Göttinger Jahre
## 2.1 Seine Studentenjahre

Während der Jahre, in denen Grandjot an der Universität von Göttingen studierte und lehrte, waren nicht nur die Mathematik in Deutschland auf der Höhe der Zeit, sondern das Göttingen Grandjots war der mathematische Mittelpunkt des Universums. Obwohl im 19. Jahrhundert in Göttingen bereits die berühmten Mathematiker Gauss, Dirichlet und Riemann wirkten, so war der große Organisator der Physik und der Mathematik jedoch Felix Klein (1849-1925). Ausgestattet mit einem großen Unternehmungsgeist und einer unglaublichen administrativen Kraft (Aufbau von Instituten, Seminaren und Bibliotheken und einer Erhöhung der Studentenzahl) entwickelte er eine Atmosphäre, die Göttingen zum „Mekka der Mathematik" machte. Hervorzuheben sind D. Hilbert, H. Minkowski, E. Landau und C. Runge, ebenso der Astronom K. Schwarzschild und die Physiker L. Prandtl, P. Debye und E. Wiechert, die für die goldene Epoche den Boden bereiteten, die am Ende der dritten Dekade des 20. Jahrhunderts ihren Höhepunkt erreichte. Studenten aus allen Ecken der Welt wetteiferten darum, in Göttingen studieren zu können. Grandjot kam 1919 nach Göttingen und alles weist darauf hin, dass er die sich dort bietenden Ressourcen maximal nutzte: so das Seminar der Physikalischen Mathematik, den Lehrstuhl der Angewandten Mathematik von C. Runge (dem ersten der angewandten Richtung in Deutschland), dem Institut der Angewandten Mathematik und Mechanik und dem Institut für Statistische Mathematik. Später beeinflusste diese interdisziplinäre Göttinger Arbeitsweise dann auch die Aktivitäten und Lehrtätigkeit Grandjots in Chile.

Die Zahlentheorie ist ohne Zweifel das, was Grandjot am meisten anzog, besonders die analytische Zahlentheorie, also das Studium der Eigenschaften der natürlichen Zahlen mit den Werkzeugen der Analysis. Diese Disziplin, von Dirichlet initiiert, hatte sich bereits in die eleganteste aber auch komplexeste Disziplin der Mathematik verwandelt. In Göttingen der 20er Jahre war ihr Spezialist par excellence der Zahlentheoriker Edmund Landau, Schüler von

Frobenius, der 1909 in Göttingen als Nachfolger von Minkowski begann und dessen grundlegendes Interesse der analytischen Zahlentheorie galt, insbesondere der Verteilung der Primzahlen. Sein „Handbuch der Lehre von der Verteilung der Primzahlen" (1909) kann als die erste systematische Darstellung der analytischen Zahlentheorie angesehen werden. Edmund Landau war der Lehrer einer ganzen Generation, nicht nur durch seine Texte und die grundlegenden Kenntnisse auf seinem Gebiet, sondern auch durch seine Schüler, unter denen z.B. P. Bernays, G. Doetsch, H. A. Heilbronn, D. Jackson, E. Kamke, A. J. Kempner, L. Neder, A. Ostrowski, W. Rogosinski, W. Schmeidler, C. L. Siegel, A. Walfisz und K. Grandjot aufragten.

1922 legte Grandjot bei seinem Doktorvater Landau die Dissertation mit dem Titel „Über das absolute Konvergenzproblem der Dirichletschen Reihen" vor. Diese von Dirichlet eingeführten Reihen erlauben die Analyse zahlentheoretischer Funktionen, besonders ihrer multiplikativen Eigenschaften.[12] Ab diesem Zeitpunkt war Grandjot Landaus Assistent und erwuchs daraus eine fruchtbare Zusammenarbeit, die fast ein Jahrzehnt andauerte.

## 2.2 Professorentätigkeit Grandjots in Göttingen

1925 wurde Grandjot in Göttingen Privatdozent[13] mit seiner Habilitationsschrift[14] „Untersuchungen über Dirichletsche Reihen" in der er Ergebnisse über Wachstumsverhalten und Nullstellen der Dirichlet-Reihen behandelt und untersucht. In den folgenden Jahren veröffentlicht er Arbeiten auf dem Gebiet der analytischen Zahlentheorie, über ganze Funktionen, Dirichlet-Reihen und trigonometrische Reihen (siehe Bibliographie).
  Die Zusammenarbeit zwischen Grandjot und Landau war sehr erfolgreich. Neben einigen Arbeiten, wie z.B. einem gemeinsamen Artikel mit Jarnik und Littlewood, arbeitete der Assistent Grandjot weit über seine Pflichten hinaus. Ganz ohne Zweifel ist eines der wichtigsten Ergebnisse dieser Kooperation das Buch „Vorlesungen über Zahlentheorie" (Leipzig, 1927, 3

---

[12] Eine gewöhnliche Dirichlet-Reihe hat die Form $\sum_{n=1}^{\infty} a_n n^{-s}$ . Eine allgemeine Dirichlet-Reihe ist eine Reihe der Form $\sum_{n=1}^{\infty} a_n e^{-\lambda_n s}$ , wobei $a_n$ komplexe Zahlen sind und $\lambda_1 < \lambda_2 < \cdots \to \infty$ .

[13] Weitere bekannte Privatdozenten in Göttingen waren A. Sommerfeld, E. Zermelo, O. Blumenthal, C. Carathéodory, E. Hecke und R. Courant.

[14] Eingereicht für die Erlangung des Titels des Außerordentlichen Professors.

Bände, 1009 Seiten), eines der einflussreichsten der Zahlentheorie überhaupt. Hardy schrieb diesbezüglich über Landau:

"This remarkable work is complete in itself; he does not assume (as he had done in the *Handbuch*) even a little knowledge of number theory or algebra. It stretches from the very beginning to the limits of knowledge, in 1927, of the "additive", "analytic", and "geometric" theories. [...] In spite of this enormous programme, Landau never deviates an inch from his ideal of absolute completeness. [...] The "Vorlesungen" is not only Landau's finest book but also, in spite of the great difficulty and complexity of some of the subject matter, the most agreeably written. The style here is the rather informal style of his lectures, which he was persuaded by his friends to leave unchanged."[15]

Im Vorwort ist Landau überaus großzügig mit den Dankesworten an seinen Assistenten Grandjot:

"Mein Dank gebührt zunächst den Verfassern der schönen Arbeiten (namentlich aus den letzten Dezennien), deren Früchte ich ernten konnte. Vor allem aber meinem langjährigen Assistenten und jetzigen Kollegen, Privatdozenten Dr. K. Grandjot, der als genauer Kenner des Gesamtgebietes mir bereits während der Vorbereitungen zu meinen Vorlesungen außerordentliche Hilfe geleistet und dann das Manuskript durchgeprüft hat. Bei den Korrekturen erfreute ich mich außer seiner Hilfe noch der Mitarbeit eines hervorragenden Forschers im Gebiete der analytischen und geometrischen Zahlentheorie, meines Schülers Dr. A. Walfisz in Wiesbaden."

Gegenüber B. Chuaqui gestand Grandjot später mit bescheidenem Stolz, dass er mindestens ein Drittel dieses berühmten Buches geschrieben habe.[16] Dass trotzdem Grandjots Name nur bei den Danksagungen erschien, muss uns heute nicht besonders überraschen, denn die Regeln der Assistenzzeiten der damaligen Epoche waren folgendermaßen: Alles, was die Assistenten produzierten, musste unter dem Namen des Chefs veröffentlicht werden.

## 2.3 Der „Einwand von Grandjot"

Unter seinen Professoren erwähnte Grandjot auch David Hilbert. Auch wenn uns leider keine Einzelheiten der Beziehung zu Hilbert bekannt sind, so ist es dennoch nicht schwierig bei seinen späteren Aktivitäten dessen Einfluss zu entdecken, denn ständig und offensichtlich zeigen sich das große Interesse und die tiefe Beherrschung Grandjots an den formalen Systemen und der Aufstellung von Axiomen.

Als zum Ende des 19. Jahrhunderts die Bemühung ihren Anfang nahm, der Mathematik eine sichere Grundlage zu geben, versicherte Kronecker: „Die ganzen Zahlen hat der liebe Gott gemacht, alles andere ist Menschenwerk!" Manche behaupteten jedoch, dass man auch die natürlichen und die ganzen Zahlen aus elementaren Bausteinen aufbauen könne. So stützte sich Richard Dedekind auf Arbeiten von Grassmann und Frege und veröffentlichte 1888 seine berühmte Abhandlung „Was sind und was sollen die Zahlen", in der er eine „algebraische Charakterisierung" der natürlichen Zahlen darstellte, die auf zwei primitiven Bausteinen aufbaute, der 1 und dem Begriff des *Nachfolgers*. Es war schließlich Giuseppe Peano, der diese Ideen fortführte und das Axiomensystem der natürlichen Zahlen[17] allgemein bekannt machte, in dem er ihnen die elegante und verständliche Form gab, wie wir sie heute kennen:

1. 1 ist eine natürliche Zahl.
2. Für jedes x gibt es genau eine natürliche Zahl als Nachfolger von x, die wir als x' bezeichnen. So ist, wenn x = y auch x' = y'.
3. Für jedes x wird angenommen, dass x' ≠ 1. Das bedeutet, es gibt keine Zahl, deren Nachfolger 1 ist.
4. Wenn x' = y', dann ist auch x = y. Das bedeutet, für jede Zahl gibt es entweder keinen Nachfolger oder der Nachfolger ist eindeutig.
5. (Induktionsaxiom) M sei die Menge der natürlichen Zahlen mit den folgenden Eigenschaften:
   a. 1 gehört zu M.
   b. Wenn x zu M gehört, dann gehört auch x' zu M.
   Dann enthält M alle natürlichen Zahlen.

Einige Jahrzehnte später erklärte Hilbert, das Ziel sei, die Analysis von den natürlichen Zahlen abzuleiten. Landau, ein Kollege von Hilbert in Göttingen, schrieb nun 1930 seine „Grundlagen der Analysis" mit dem ausdrücklichen Vorsatz für einen Text zu sorgen, der die Analysis anhand der Axiome von Peano entwickelt, ein Vorhaben, dem sich bisher noch niemand gestellt hatte. Wörtlich formulierte dies Landau wie folgt: „Nun gibt es in der ganzen Literatur kein Lehrbuch, das sich das bescheidene Ziel setzt, nur das Rechnen mit Zahlen im obigen Sinne zu begründen." und gab an, dass er mit seinem Buch anträte, diese Arbeit zu erledigen, da es sonst noch niemand getan hätte. Bevor sie zu

---

einem Buch wurden, gab Landau seine Aufzeichnungen seinem Assistenten, der darüber eine Vorlesung hielt. Zum Schluss des Kurses erhielt Landau die Aufzeichnungen zurück mit ausführlichen schriftlichen Anmerkungen, die im Wesentlichen aussagten, dass es nicht möglich sei, nur mit den Axiomen von Peano die gesamte Analysis abzuleiten und fügte auch sogleich die Axiome bei, die das Problem lösten. Dies wurde als „Einwand von Grandjot" berühmt. Das Problem erklärte Landau mit seiner ihm eigenen Deutlichkeit:

„Wenn ich etwa in einer Vorlesung über Zahlentheorie irgend einen Satz über natürliche Zahlen so beweise, dass ich erst die Richtigkeit für 1 und dann aus der Richtigkeit für x die für x + 1 beweise, so pflegt gelegentlich ein Zuhörer den Einwand zu erheben, ich hätte die Behauptung ja gar nicht vorher für x bewiesen. Der Einwand ist unberechtigt, aber verzeihlich; der Student hatte eben nie vom Induktionsaxiom gehört.

Grandjots Einwand klingt ähnlich; mit dem Unterschiede, dass er berechtigt war, so dass ich ihn auch verzeihen musste. Auf Grund seiner fünf Axiome definierte Peano die Summe x + y bei festem x für alle y folgendermaßen:
x + 1 = x',
x + y' = (x + y)',
und er und Nachfolger meinen damit: x + y ist allgemein definiert; denn die Menge der y, für die es definiert ist, enthält 1 und mit y auch y'. Aber man hat ja x + y gar nicht definiert."

Grandjot löste das Problem, indem er zusätzliche Axiome hinzufügte. Landau bevorzugte letztendlich einen Vorschlag des ungarischen Logikers Kalmár.[18] Der berühmte Einwand von Grandjot wird immer wieder kommentiert. Eine seiner relevantesten Folgen ist die scharfsinnige Beobachtung über die Grenzen der einfachen Induktion und die Rolle, die Definitionen wie $\sum_i x_i$ und $\prod_i x_i$ im Aufbau der Analysis spielen.[19]
   Um diesen Abschnitt abzurunden, möchten wir noch bemerken, dass Grandjot während seiner Zeit in Göttingen viele Kontakte knüpfte, die er auch noch in Chile pflegte: wir möchten da G. Birkhoff, G. H. Hardy und F. Hausdorff unter vielen anderen hervorheben. Diese Mathematiker waren später die Verbindungspunkte, durch die Grandjots Schüler von Chile aus an die weltweit wichtigsten Zentren gelangen konnten, um Mathematik zu studieren.

## 3 Die Ankunft in Chile

### 3.1 Das wissenschaftliche Umfeld in Chile bei seiner Ankunft

In den 20er Jahren befand sich Chile in großen Turbulenzen politischer, sozialer und universitärer Art, die alle Institutionen, auch den universitären Lehrkörper durchdrangen und überall spürbar waren. Rektoren genau wie Minister der Regierung folgten in einer solchen Geschwindigkeit aufeinander, wie niemals zuvor. Zwischen 1926 und 1930 hatte die Universität von Chile nicht weniger als 5 Rektoren: Amunátegui, Matte, Charlín, Martner und Quezada, die „mit allen zur Verfügung stehenden Mitteln versuchten, das studentische Durcheinander zu kontrollieren".[20] Der Aufruhr durchdrang alle Bereiche; ein großes Interesse an der Entwicklung der nationalen Wissenschaften war geweckt worden. Die Steigerung und eine breitere Fächerung der nationalen Industrie, die Suche nach neuen politischen Formen der Bürgervertretung und die pädagogischen Reformen, die 1931 ihren Höhepunkt mit der Aufstellung der neuen Statuten der Universität von Chile erreichten, trugen wahrscheinlich dazu bei.

Rektor Martner war es wichtig, weit über eine bloße Berufsausbildung hinauszugehen und sich den wissenschaftlichen Entwicklungen zu nähern und er sagte 1928: „Auch wenn viele es gerne so beibehalten wollen, ist der kulturelle Auftrag einer Universität im Wesentlichen nicht die Versorgung mit Wissen, das schon lange bekannt und erforscht ist oder das schon Bekannte zu veröffentlichen, *sondern als Quelle der Forschung zu dienen und ein Hebelpunkt für die Entwicklung der Wissenschaften zu sein.*" (Hervorhebung mit Kursivschrift durch die Autoren). Weiter führte er aus: „Es ist notwendig, Seminare, Forschungslabore und spezialisierte Bibliotheken einzurichten, so dass jeder Lehrstuhl sein eigenes Seminar, Labor oder eigene Bibliothek als unterstützendes Arbeitsmittel besitzt und damit die höheren Studien erfolgreich absolvieren kann."[21]

---

[18] Die zusätzlichen Axiome von Grandjot sind uns nicht bekannt. Landau kommentierte den Vorfall – er nannte ihn „Abenteuer" – ausführlich im Vorwort seiner „Grundlagen der Analysis".
[19] Der Einwand von Grandjot hat noch viel mehr an Logik und Mathematik zu bieten. Daher empfehlen wir dem interessierten Leser den Artikel des Prof. Pi Calleja, „La objeción de Grandjot", Mathematica Notae, 10. Jahrgang, 1940, S. 143 – 151.

[20] Rolando Mellafe, „Reseña histórica del Instituto Pedagógico", Universidad Metropolitana de Ciencias de la Educación, 1988, S. 13.
[21] Textos Universitarios, in Daniel Martner U., *Obras Escogidas (Ausgewählte Werke)*, Edición del Centro de Estudios Político Latinoamericanos Simón Bolívar, 1992

Innerhalb der Fakultäten führte diese Bewegung zu ersten Ergebnissen. 1928 konnte die Fakultät der physikalischen Wissenschaften den französischen Physiker Paul Langevin und den italienischen Logiker Federico Enriques begrüßen, die Konferenzen über ihre Arbeitsgebiete abhielten und später feierlich als Ehrenmitglieder in die Fakultät aufgenommen wurden.[22] Berühmt waren damals auch die Konferenzen des Ingenieurs Ramón Salas Edwards über die Relativitätstheorie im Jahr 1929. All das zeigt die Bemühungen der Fakultät, dem Fachwissen eine Basis auf der Höhe des neuesten wissenschaftlichen Standes zu geben.

Als weitere Auffälligkeit, die die Bewegungen in der Entwicklung der Wissenschaft in dem damaligen Jahrzehnt hervorbrachte, ist die Eröffnung des Institutes der Wissenschaften Chiles (Instituto de Ciencias de Chile), das dazu bestimmt war, Forschungen und reine wissenschaftliche Studien zu koordinieren und zu fördern, ohne zunächst auf einen bestimmten Zweck gerichtet zu sein, um damit die Kultur sowohl zu bewahren als auch zu erhöhen und die größten nationalen Probleme lösen zu können."[23] Dieses Institut wurde durch 3 Akademien gebildet: Akademie der Ökonomischen und Sozialen Wissenschaften, Akademie der Mathematik und Naturwissenschaften und Akademie der Geschichte, Philosophie und Philologie. In der Akademie der Mathematik und Naturwissenschaften findet man als Gründerpersonen die Namen Enrique Froemel, Ricardo Poenisch und Carlos Grandjot, also drei wichtige Persönlichkeiten im Prozess der Entwicklung der Mathematik in Chile: Poenisch und Froemel waren dabei besonders im schulischen und universitären Mathematik tätig und Grandjot in der wissenschaftlichen Mathematik und in der Aus- und Weiterbildung von Schülern nach seinen Prinzipien; alle drei zusammen in der angewandten Mathematik der Ingenieure.

Die Stimmung, die am Ende der 20er Jahre in Chile herrschte, wird rückblickend gut zusammengefasst in einer akademischen Abhandlung von Juvenal Hernández, einem ehemaligen Rektor der Universität: „Als Folge dieser Ereignisse (Turbulenzen) die die globalen Erschütterungen (1. Weltkrieg, Industrielle Revolution) widerspiegelten, begann unser Land seine Unterentwicklung abzuschütteln. Es begann zu knospen und zu kämpfen; ein gewaltiger Schub von Schöpfung und Erneuerung setzte ein."[24] „Als ich 1933 das Amt des Rektors übernahm", erinnert sich Hernández weiter, „war die Universität berufsbildend und akademisch (…), es war nötig, sie zu entwickeln. Ich habe einige Initiativen angestoßen, die einen Anreiz für die reine Forschung und die Anwendung der wissenschaftlichen Erkenntnisse gegeben haben (…), und daraus entwickelten sich zum ersten Mal Institute, Seminare, Werkstätten und Laboratorien." Tatsächlich war diese Amtszeit Hernández' (1933–1953) fruchtbar bei der Hervorbringung von Körperschaften, die die wissenschaftliche Forschung förderten. So entstanden trotz der Krise in den 30er Jahren und trotz der negativen Auswirkungen des 2. Weltkrieges auf die peripheren Länder wie Chile dreißig Institute der verschiedensten Forschungsgebiete. Bevor diese Umwälzungen stattfanden, befand sich die organisierte Forschung in einem Dämmerzustand, was aber nicht bedeutete, dass es nicht auch eine wissenschaftliche Produktion mit guter Qualität gab. Als Beispiel moderner Forschung vor den 30er Jahren kann das Institut der Physiologie der Universität in Concepción gelten und das Bakteriologische Institut der Medizinschule der Universität von Chile, die beide, wie A. Meyer sagte, „keinen Vergleich mit ähnlichen Institutionen, die es auf der ganzen Welt verstreut gibt, scheuen müssen."[25] Es waren moderne Forschungseinrichtungen, aber es dominierte dort das Aspekt des Lehrens als herausragendes Ziel.

## 3.2 Grandjots Vorgänger

Die Politik Chiles, ausländische Gelehrte zu verpflichten, um Wissenschaft ins Land zu holen, war durchgehende Praxis seit der Unabhängigkeit, als das chilenische Schulsystem nach dem Vorbild des französischen strukturiert wurde. Die ersten Texte der schulischen und universitären Mathematik waren wörtliche Übersetzungen der französischen Vorlagen und auch die späteren von nationalen Autoren folgten eng der französischen Didaktik.[26] Die französische Schule, die die Lehre der Mathematik in Chile seit den ersten Jahrzehnten des 19. Jahrhunderts bestimmte, machte ab 1889 der deutschen Schule Platz. In diesem

Jahr kam die bereits erwähnte Gruppe deutscher Professoren in Chile an, die von der Regierung verpflichtet waren, das Pädagogische Institut der Universität Chiles mit Leben zu erfüllen. Dieses Institut war per Dekret vom 29.4. 1889 gegründet worden mit dem Ziel, Mathematiklehrer der höheren Schulen auszubilden.

Die Anstellung Grandjots und anderer Professoren[27] im Jahr 1929 folgte aus dem Wunsch der chilenischen Regierung, die erste Generation derjenigen zu verstärken, die sich in den 4 vorherigen Jahrzehnten intensiver Entwicklungsaktivitäten abgearbeitet hatten. Die Lehrmeister von 1889, so Enrique Molina, „waren praktisch ohne Ausnahme wirkliche Männer einer fleißigen Wissenschaft, das begannen und völlig ihren Studien hingegeben."[28] Einige von ihnen, wie Federico Johow, Rodolfo Lenz und Ricardo Poenisch, machten Chile zu ihrer zweiten Heimat. Dasselbe taten später Carlos Grandjot und Ferdinand Oberhauser.[29]

Sowohl die schulische als auch die universitäre Mathematik waren bei der Ankunft Grandjots bereits gut organisiert und hatten dank der Arbeit von Tafelmacher und Poenisch einen guten Ruf. Seit Beginn hatten sie die Verantwortung für die Ausarbeitung der Programme und die Abfassung der mathematischen Lehrtexte sowohl für die unteren als auch für die höheren Schulen.[30] Augusto Tafelmacher bildete die erste Generation von Mathematiklehrern am Pädagogischen Institut aus. Poenisch folgte ihm auf diesem Lehrstuhl 1908.[31] Ende der 20er Jahre verbreiteten dessen Schüler die Mathematik im ganzen Land durch ihren Unterricht in Lyzeen, staatlichen Schulen, Privatschulen und weiteren Lehranstalten. Dabei verließ die Mathematik die Grenzen einer nützlichen Wissenschaft (*ciencia "util"*), in der sie lange gehalten wurde und wurde zu einer kulturellen und autonomen Disziplin, ein Unternehmen, das in Frankreich, Deutschland und anderen europäischen Ländern schon seit dem frühen 19. Jahrhundert realisiert worden war.

Als Erfolg der Kraftanstrengungen der Professoren von 1889, gab es nun auch bald landeseigene Lehrer für den Unterricht in Höheren Schulen. Die besten Schüler von Poenisch folgten ihm sowohl am Pädagogischen Institut, als auch in der Ingenieur- und der Militärschule.[32] Andere dienten in den staatlichen Schulen, der Schule für Schifffahrt, der Luftfahrtschule, in der Schule für Kunst und Gewerbe, in den neu entstehenden Privatuniversitäten und in anderen Schulen der Universität von Chile wie z.B. der Schule für Landwirtschaft und der für Architektur.

Auch die Physik, das zweite Spezialgebiet von Grandjot, verstärkte gemeinsam mit der Mathematik ihre schulische Lehre. Ziegler und Gostling, beide Lehrstuhinhaber am Pädagogischen Institut, achteten sehr darauf, dass diese Disziplin nach allen methodischen Regeln einer experimentellen Wissenschaft unterrichtet wurde. Man erzählte sich, dass immer, wenn in irgendeiner Provinz ein Lyzeum neu errichtet wurde, Ziegler beim Schulministerium vorsprach und die Einrichtung eines Physiklaboratoriums einforderte.[33] So groß war seine Besorgnis um die Weiterentwicklung dieser Wissenschaft, auch noch in den letzten Jahren seiner unterrichtlichen Tätigkeit. Die Texte *Experimentalphysik* von Ziegler und Gostling, zu Beginn des Jahrhunderts geschrieben, erreichten 1952 ihre 13. Auflage.

Nach Meinung von Carlos Videla, einem Schüler von Poenisch und einer seiner Nachfolger,

---

[27] Verpflichtet von der chilenischen Regierung kamen 1929 außer Grandjot: Ferdinand Oberhauser für Chemie; Guillermo Goetsch für Biologie, Adolph Meyer für Philosophie; Woldemar Voigt für Praktische Pädagogik; Peter Petersen als Fachleiter für Höhere Schulbildung, außerdem der nordamerikanische Ovied Hundley als Laborchef der Metallurgie. Boletín Universitario U. Chile, 1929.

[28] E. Molina: El primer curso del Instituto Pedagógico (Der erste Kurs des Pädagogischen Institutes), zum 75. Jahrestag seiner Gründung, Universität von Chilöe, 1964. Molina war Schüler dieses Kurses.

[29] Auch wenn Poenisch nicht Gründungsprofessor del Pädagogischen Institutes war, gehörte er trotzdem zu dessen erster Generation. Er kam 1889 nach Chile und vor dem Pädagogischen Institut erfreuten sich die Ingenieursschule der Universität von Chile, das Nationale Institut und das Lyceum von Rancagua seiner Dienste.

[30] Für die Grundschulen erarbeiten sie 2 Bände für Geometrie, 2 für Algebra, 1 für Trigonometrie und 1 für Stereometrie, die mit einigen Abänderungen 70 Jahre in Chile verwendet wurden. Für die Höheren Schulen schrieb Tafelmacher die Abhandlung der sphärischen Trigonometrie, Elemente der Analytischen Geometrie und Elemente der Höheren Algebra (Tratado de Trigonometría esférica, Elementos de Geometría Analítica, Elementos de Algebra Superior). Die Arbeiten von Poenisch sind: Ein Text der Analytischen Geometrie (Geometría Analítica) und ein anderer über Analysis (Análisis), der die höhere Algebra einschließt. Er publizierte sie unter dem Titel Einführung in die Höhere Mathematik (Introducción a las Matemáticas Superiores) und sie dienten viele Jahrzehnte als Studientexte für die Schüler des Päd-

agogischen Institutes, der Militärschule und der Ingenieurschule der Universität von Chile.

[31] Der Professor, der 1889 eigentlich für Mathematik verpflichtet worden war, war Dr. Ricardo von Lilienthal, der allerdings nur sehr kurz in Chile blieb. Ab 1890 wurde er durch Tafelmacher ersetzt.

[32] Hier sind einige Namen: Carlos Videla, Enrique Froemel, Abraham Pérez, Oscar Marín, Jenaro Moreno, Domingo Almendras, Federico Rutland, Agustín Rivera. Alle waren bis in die 50er Jahre Kollegen von Grandjot. Mit einer kurzen Zeitverschiebung muss man auch noch Guacolda Antoine und César Abuadad hinzufügen.

[33] Erinnerungen von Prof. Raquel Martinolli, Schülerin von Ziegler und seine Assistentin im Jahr 1929 (Decreto 11. Februar 1929, Boletín Universitario).

wussten diese Lehrmeister „ihre Arbeit an die Bedürfnisse im Lande anzupassen und die Möglichkeiten auszunutzen, die der kulturelle Grad, der in dieser Epoche erreicht wurde, zuließ. Ihre Publikationen[34] verrieten Enthusiasmus und Begabung sowohl für die Forschung als auch für die Information und die umfassende Beherrschung der Mathematik.“[35] Es besteht kein Zweifel, dass ihr Auftrag, Lehrer für die Schulen auszubilden, ein voller Erfolg war. Auf dem Gebiet der Forschung war der Erfolg nicht ganz so groß, trotz Enthusiasmus und Begabung. Die systematische Ausbildung von Nachfolgern war eine Aufgabe, die ein neuen „Import" von gelehrten Ausländern bzw. dem Erwachen von nationalen Talenten vorbehalten war.[36]

## 4 Grandjots neue Heimat
### 4.1 Die ersten Jahre: 1929-1945 Axiomensysteme und Abstrakte Algebra

Bei der Ankunft in Chile begriff Grandjot sehr schnell den Geist der Veränderung und der Reform, von dem die Gesellschaft Chiles ergriffen war. Dabei zeigte sich besonders die akademische Schicht daran interessiert, ihm in der Universität die nötigen Voraussetzungen zu schaffen, um eine moderne Forschung zu ermöglichen.

In seiner kurzen Biografie erzählt Grandjot, dass er zunächst im Pädagogischen Institut der Universität von Chile mit Kursen zur Höheren und Elementaren Mathematik, zur Philosophie und zur Physik begann. Parallel dazu bot er Seminare an, in denen entweder mathematische Artikel gelesen und diskutiert oder Themen, die mit der Lehrtätigkeit zusammenhingen, vertieft wurden. An diesen Zusammenkünften, die immer an den Samstagen stattfanden, nahmen auch zahlreiche der damaligen jungen Professoren, Schüler von Poenisch teil: Videla, Marín, Almendras, Pérez und andere.[37]

Die Mathematik der damaligen Zeit in Chile verharrte bei einem „klassischen" Stil. Weder Mengenlehre, Axiomensysteme, noch moderne oder abstrakte Algebra wurden unterrichtet. In den europäischen Ländern waren diese Themen bereits Bestandteil der Unterrichtsprogramme. Für das vierjährige Studium am Pädagogischen Institut wurden folgende Kurse angeboten: elementare Mathematik, analytische Geometrie, sphärische Trigonometrie, klassische und höhere Algebra, und Differenzial- und Integralrechnung. Das Rechnen wurde bis zum Ende des 19. Jahrhunderts mit Hilfe der Texte von Francoeur unterrichtet. Geschrieben wurden diese Texte zu Beginn des Jahrhunderts im Stil von Newton oder Leibnitz, das heißt, ohne die Konzepte von Funktion, Grenzwert, Stetigkeit oder der Ableitung von Funktionen nach Cauchy; erst Poenisch führte diese Konzepte in seinen Analysis-Kursen ein. Ab 1930 waren diese Ideen in der chilenischen Mathematik weit verbreitet und wurden von den Professoren beherzt angewandt. Allerdings blieben den unendlichen Reihen, die für die Infinitesimalrechnung unentbehrlich sind, eine strenge Anwendung verwehrt. In einem Artikel, den Grandjot über dieses Thema in der *Revista de Matemáticas y Física*, die in Chile weitverbreitet war, publizierte, schrieb er: „Nachdem ich einige Unterrichtstexte gründlich untersucht habe, muss ich feststellen, dass es dort keine Vorschläge gibt, wie mit dieser Materie umzugehen sei. [...] Einige Autoren streben so sehr danach, sich auf das Niveau der Schüler herunter zu begeben, dass sie die Genauigkeit der Mathematik dabei opfern, bei anderen Autoren wiederum bemerkt man ganz klar, dass sie die Materie gar nicht verstanden haben." In seinem Artikel bemüht sich Grandjot darum, korrekte Definitionen vorzugeben und grundlegende Sätze für die Theorie der Reihen aufzustellen. Dabei ging er streng vor, aber so, dass die Theorie von allen verstanden wurde, „bis hinunter zu den Anfängern".

Zusätzlich zu diesen Themen aus der Analysis führte Grandjot im Pädagogischen Institut 1936[38] Kurse über Differentialgeometrie und über die Grundlagen der Mathematik bzw. der Axiomatisierung ein, die bis zur 1945er-Reform der Philosophischen Fakultät, zu der das Pädagogische Institut gehörte, stattfanden. Der Kurs der Differentialgeometrie führte zur Vektor-Analysis und der Grundlagen-Kurs blieb ein Wahlkurs in der Lehrerausbildung der Mathematiker und Physiker. Diese Reform verlängerte die Ausbildung auf 9 Semester.

---

[34] Siehe: *Anales de la Universidad de Chile*: 1892, 1893, 1894, 1897, 1901 und 1905.
[35] Carlos Videla, *Contribución de la Facultad de Filosofía y Humanidades a la enseñanza de las Matemáticas en Chile*, 1944, in *El Centenario de la Universidad de Chile*.
[36] Die schöpferische Unruhe war trotzdem vorhanden. In den ersten Jahrzehnten des 20. Jahrhunderts registrierten die Annalen der Universität von Chile Publikationen über mathematische Themen, die von einheimischen Autoren verfasst wurden, so über Nepersche Logarithmen, Gleichungen 2. Grades oder mehrfache Integrale. Siehe *Anales* 1925, 1927 und 1930, Band VIII.
[37] Erinnerungen der Professorin Guacolda Antoine, aus Gesprächen mit den Autoren. Antoine, die 1928 habilitierte, nahm an diesen Abenden als Gast teil.

---

Sie erinnerte sich besonders an die Sitzungen, die partielle Differentialgleichungen thematisierten.
[38] Siehe sein Arbeitstagebuch in der Universität von Chile.

Der Grundlagen-Kurs war der neuartigste Kurs in unserem damaligen wissenschaftlichen Umfeld. Berücksichtigt man die entsprechenden Verhältnisse, dann war er dem von Hilbert in Göttingen über die Grundlagen der Geometrie ähnlich, der als Antwort auf den Intuitionismus konzipiert worden war. In einer auf Axiomen aufbauenden Theorie interessiert nicht die eigene Natur der Dinge; das, was wichtig ist, sind die Relationen zwischen ihnen. Dieser Kurs lenkte die Aufmerksamkeit der Studenten auf die Definitionen der Begriffe und der Axiome mit denen man sich normalerweise im Unterricht der Mathematik auseinandersetzt. In der damaligen Zeit war es üblich, in den Schul- und sogar in den Lehrbüchern Definitionen wie die folgende zu lesen:

„Eine Zahl ist das Ergebnis des Vergleichs einer Menge oder Größe mit einer Einheit." „Eine Gerade ist der kürzeste Abstand zwischen zwei Punkten." Für die Verständlichkeit dieser Definitionen hätten zunächst einmal die Ausdrücke „Menge", „Einheit" und „Abstand" definiert werden müssen, Begriffe, die ihrerseits nur definiert werden können mit Aussagen über „Zahl" im ersten Fall und „Gerade" im zweiten Fall.

Diese Definitionen sind also jeweils Teil eines Teufelskreises. Um diesen Kreis zu durchbrechen, werden bei der Axiomatisierung bestimmte Termini als grundlegend bestimmt, die ohne Definition angenommen werden und von diesen ausgehend definiert man alle weiteren Begriffe, die nötig sind. Der zweite Schritt bei der Aufstellung eines Axiomensystems besteht in der Darlegung einer Menge von Lehrsätzen, die man Axiome nennt und die man ohne Beweise akzeptiert. Von diesen leiten sich weitere Lehrsätze ab. Im dritten Schritt ist es notwendig, Regeln aufzustellen, die es ermöglichen von den Axiomen neue Ableitungen zu bilden, die man Theoreme nennt. Axiomensysteme müssen drei Eigenschaften erfüllen: Widerspruchsfreiheit, Unabhängigkeit und Vollständigkeit, die hier nicht ausführlich behandelt werden sollen, auf die Grandjot aber großen Wert legte und die er sorgfältig in einem ad-hoc konstruiertem Modell bewies.[39] Erinnern wir uns besonders an den Kurs der Grundlagen der Geometrie (1951), bei dem einer der grundlegenden Begriffe die Beziehung „zwischen" ist. Eine der Übungen, die Grandjot in diesem Fall vorschlug war folgende: „Beweist, dass die Summe der Winkel in einem Dreieck kleiner als 180° ist, indem ihr diese Lehrsätze berücksichtigt (die er aufzählte)."

---

[39] Obwohl Eigenschaften wie die Unabhängigkeit eines Axiomensystems wichtig sind, konzentriert sich die moderne Logik auf zwei Bedingungen: Korrektheit, das System darf also keine unkorrekten Lehrsätze ableiten und Vollständigkeit, was bedeutet, dass das System alle Lehrsätze vollständig ableiten muss.

Aber er ermunterte nicht nur seine Schüler zu Übungen in der reinen Mathematik, sondern er lenkte ihre Aufmerksamkeit auch auf damit zusammenhängende Probleme. So bestand eine weitere seiner Übungen darin, dass er Schulbücher darauf hin untersuchen ließ, wie streng sich an die mathematischen Prinzipien der Axiomatisierung hielten. Ergebnis genau dieser Analyse war die „Entdeckung" der oben beschriebenen Teufelskreise bei der Definition von „Zahl" und „Gerade". Ebenfalls Ergebnisse dieses Kurses waren auch einige Abschlussarbeiten über Zahlensysteme und nichteuklidische Geometrie von Studenten des Pädagogischen Institutes.

Die Prinzipien der Axiomatisierung sind ein wichtiges Werkzeug auf allen Feldern des Wissens und auch im täglichen Leben: „Wenn Sie mit mir diskutieren wollen", so sagte schon Voltaire, „dann definieren Sie alle Begriffe, die Sie verwenden." Daher hatte der Kurs auch von Studenten der Philosophie, des Rechts und anderer Disziplinen einen großen Zulauf, denn dort ist die deduktive Beweisführung grundlegend.

In einem anderen Abschnitt seiner Lehrtätigkeit, 1933, begann er einen Kurs „Ergänzungen der höheren, reinen und angewandten Mathematik" in der Ingenieurschule der katholischen Universität von Chile. Von 1945 bis 1963 bot er diesen Kurs auch in der Ingenieurschule der Universität von Chile an. Die Vorlesungsinhalte veröffentlichte er 1950 im Editorial Universitaria in 2 Bänden zu jeweils 300 Seiten. Diese unterteilen sich in 4 Abschnitte: graphische und numerische Methoden, Analysis, Differentialgleichungen sowie partielle Differentialgleichungen der mathematischen Physik. Diese Bände waren dazu bestimmt, die mathematischen Erkenntnisse der letzten Jahre zusammenzufassen. Im Vorwort stellt der Autor fest: „Wegen der Stellung zwischen der reinen und der angewandten Mathematik musste ich einen vernünftigen Kompromiss für die Strenge der Beweisführung finden." Der erste Teil behandelt das Thema auf der Grundlage von Beispielen und Grandjot zeigt dabei anhand einer Sammlung von Rechentricks seine große Meisterschaft im Kopfrechnen. Der Rest ist ausgeglichen im Gebrauch der Intuition und einer „strengen Beweisführung". Es ist durchaus möglich, dass die Ideen zu diesem Kurs bis an den Lehrstuhl der angewandten Mathematik in Deutschland zurückreichen, wo Grandjot bei C. Runge ein herausragender Schüler war. Für die Lehre in Chile war dieser Beitrag von großem didaktischem Wert.

Sein Eifer, die Schüler an die aktuellsten Themen heranzuführen, trieb ihn dazu, eine Monographie der abstrakten Algebra zu schreiben, die er etwa 1940 in der Revista Universitaria de la Pontificia Universidad Católica publi-zierte. Wäre dieser Text dann auch in der Lehre verwendet worden, dann hätte er, unserer Meinung nach, das offizielle Studium dieser „neuen Wissenschaft" in Chile um einige Jahrzehnte beschleunigen können. Die moderne oder abstrakte Algebra ist eine Disziplin, die sich zu Beginn des 20. Jahrhunderts entwickelte, und sie ist die systematische Weiterentwicklung der Verallgemeinerung arithmetischer Operationen mit Hilfe von Symbolik und mit Axiomatisierung. Sie untersucht die verschiedenen Typen der algebraischen Strukturen, die als Verallgemeinerung (Abstraktion) von konkreten Strukturen wie Transformations-Gruppen oder die Matrix- und Idealtheorie auftreten. Die eigentlichen Ursprünge dieser Disziplin kann man 1910 an der Publikation von Steinitz' „Algebraische Theorie der Körper" festmachen. Danach, 1926, veröffentlichte L. E. Dickson seine „Modern Algebraic Theories" und 1936 erschien die klassische, wegweisende Abhandlung „Moderne Algebra" von B. L. van der Waerden. Zu dieser Zeit existierte keine Arbeit in spanischer Sprache über das Thema und die oben genannten waren schwer zugänglich. Der klassische Text, der in den Vereinigten Staaten verwendet wurde, „Modern Algebra" von Birkhoff und Mac Lane, war aus dem Jahr 1941, und die ins Spanische übersetzte Version stammt aus dem Jahr 1954. Grandjot schreibt: „Auf den folgenden Seiten werde ich versuchen, die Elemente dieser neuen mathematischen Wissenschaft aufzuzeigen. Ich werde ihre grundlegenden Konzepte darlegen; aber die riesigen Anwendungsmöglichkeiten kann ich nur anreißen." Grandjots Text überrascht mit seiner Modernität und seinen Visionen. Er begann mit einem Abschnitt über Mengenlehre, einem Thema, das erst Jahrzehnte später im schulischen Unterricht Chiles berücksichtigt wurde. Er fuhr fort mit einer axiomatischen Abhandlung der Theorie der Ringe und Körper, einer ausgesprochen modernen Problemstellung. Die Form seiner Ausführungen ist sehr vorsichtig gehalten. Vor jedem neuen Thema führt er eine Begründung mit konkreten Beispielen an, denn Grandjot war sich sehr bewusst, dass er die Themen Lesern nahebringen musste, die abstrakte Beweisführungen nicht besonders gewohnt waren. Die Ergebnisse dieser Methodik waren außergewöhnlich.

Trotzdem wurde dieser Text, der dazu gedacht war, Chile mit den abstraktesten und zukunftsträchtigsten mathematischen Disziplinen bekannt zu machen, im universitären chilenischen Umfeld eisig aufgenommen, um nicht zu sagen: gar nicht beachtet. Weder das Haus der höheren Studien (casa de estudios superiores), in dem er veröffentlicht wurde, noch die Universität von Chile, die damals verantwortlich war für Ausrichtung und Kontrolle aller Studiengänge, und noch nicht einmal die Kollegen der eigenen Fakultät äußerten sich zu dieser exzellenten Arbeit der modernen Mathematik. Währenddessen sen in Argentinien und Brasilien die abstrakte Algebra in die Lehrprogramme schon Mitte der 30er Jahre aufgenommen wurde[40], gelang es Grandjot trotz größter Anstrengungen in Chile nicht, diese neue Wissenschaft im Land bekannt zu machen. Hatte man der chilenischen wissenschaftlichen Entwicklung die Tradition und den Konservatismus aufgezwungen? Oder eine Routine verordnet? Was auch immer die Antworten auf diese Fragen sein mögen, sicher ist, dass wir bis in die 70er Jahre warten mussten, bis die moderne Algebra als regulärer Zweig in die Lehrprogramme aufgenommen worden war und auch das ist nur der Hartnäckigkeit Grandjots und den Anstrengungen seiner Schüler zu verdanken. Da aber die damalige Regierung in ihren Entwicklungsplänen keinen Platz für avantgardistische Wissenschaften gelassen hatte, ihre Weiterentwicklung zu fördern und sie zu einem lang anhaltenden Wachstum zu ertüchtigen, geriet die „Algebra Abstracta" von Grandjot viele Jahre wieder in Vergessenheit.

Aber Grandjot kümmerte sich nicht nur um die Perfektionierung des Hochschulunterrichtes, sondern es war ihm auch wichtig, eine Verbindung zwischen Chile und Europa auf der Ebene des Schulunterrichtes zu schaffen. Dafür entwickelte er in Zusammenarbeit mit dem Lehrer und späterem Professor Oscar Marín einen Arithmetikkurs für die ersten Jahre der mittleren Schulbildung, von dem allerdings leider nur der erste Band publiziert wurde. Das Buch war so strukturiert, dass es das Selbststudium des Schülers stark unterstützte: Jeder Absatz begann mit anschaulichen Beispielen aus dem Lebensbereich der Schüler, ging über zur Theorie und kam zurück auf die praktischen Anwendungen, um nun umfassenden sozialen Interessen oder dem reinen Intellekt zu dienen. Die alten Methoden, die nur auf Auswendiglernen basierten, waren damit überwunden. Diese methodische Vorgehensweise lag ganz auf der Linie des Dalton-Plans bzw. des Winnetka-Systems, welches in den 30er Jahren in Mode war.[41] Auch dieses Buch geriet, wie seine *Ab-*

---

[40] Luigi Fantappiè, 1934, in Brasilien, und Sagastume Berra, 1937, in Argentinien. Siehe: Pereira da Silva, A matemática no Brasil, Curitiba, Ed. Da UFPR, 1992, S. 235. Ebenso: Sagastume Berra, Lecciones de Algebra moderna, Argentinische Republik, La Plata, 1961, Vorwort.

[41] Im Englischen Institut in Santiago de Chile, welches direkt beim Pädagogischen Institut steht und wo heute die Hauptstadtuniversität der Erziehungswissenschaften beherbergt wird, wurde der Dalton-Plan mit Erfolg seit 1929 für einige Zeit angewandt. Es ist sehr wahrscheinlich, dass das Buch von Grandjot und Marín in direktem Zusammenhang mit dem Dalton-Experiment stand.

*strakte Algebra,* für mehr als ein halbes Jahrhundert in Vergessenheit.

In dem Zeitraum zwischen seiner Ankunft in Chile und dem Ende des Zweiten Weltkrieges gab es mehrere Schicksalsschläge, die Grandjots Zukunft beeinflussten. Trotzdem war wahrscheinlich diese Zeit diejenige, die er am meisten im Kreise seiner Familie genoss. Mit seiner Frau Gertrudis erforschte er einen großen Teil des chilenischen Territoriums, immer auf der Suche nach Pflanzen und autochtonen Arten. Im Jahre 1941 unternahmen sie sogar eine Reise durch Bolivien und Hoch-Peru. Man sagte, dass beide Experten auf dem Gebiet der Botanik waren. Sie berichteten von ihren Entdeckungen auf Zusammenkünften der chilenischen Vereinigung der Naturgeschichte, in der sie Mitglieder waren und publizierten sie in den entsprechenden Fachzeitschriften. Nibaldo Bahamonde, Preisträger des Nationalpreises der Wissenschaften, erinnert sich an seine Studienzeit im Pädagogischen Institut in den 40er Jahren. Bei Ausflügen seines Chemiekurses habe „Professor Grandjot, ein ausgezeichneter Mathematiker, überaus sympathisch, mit großer Liebe zur Natur und überragender Kenntnis der Botanik" für gewöhnlich diese Exkursionen zu vielen Orten des Landes gemeinsam mit Professor Oberhauser begleitet.[42]

Grandjot besuchte Deutschland nur selten. In dem Zeitraum, den wir betrachten, besuchte er seine Heimat in den Sommern von 1931 und 1938. Dabei fand er auf seiner letzten Reise die Universitäten unter strenger Kontrolle vor. Der „Glanz und die außerordentliche Ausstrahlung", die die deutsche Mathematik zwischen 1920 und 1933 hatte, „waren auf eine brutale Weise verstümmelt worden".[43] Ohne Zweifel war das traurigste für Grandjot, dass sein Doktorvater Edmund Landau seiner Ämter enthoben wurde und im Ausland war. Diese schmerzhaften Wunden wurden später noch vergrößert durch die vom Weltkrieg erzwungene Trennung von seinen Gleichgesinnten in Deutschland. Auch in Chile gab es als Folge dieses großen Konfliktes schwierige Momente zu ertragen. Obwohl er nicht mit dem Naziregime sympathisierte, wurde er seiner Ämter im Pädagogischen Institut enthoben und nach Rengo „verbannt". Dort bekam er Typhus und wegen der prekären medizinischen Zustände in Rengo wurde er zur Behandlung in der Deutschen Klinik in Santiago aufgenommen. Um danach seine Ämter in der Universität wiederzubekommen, musste er sich gemeinsam mit anderen Anwärtern bewerben, darunter auch einigen seiner ehemaligen Schüler, was ihn besonders schmerzte. Dank seiner großen Verdienste und der Hilfe seiner früheren Kollegen konnte er seine akademischen Ämter nach den alten Konditionen wieder aufnehmen.

Und als wenn all das Genannte noch nicht schlimm genug gewesen wäre, so kam 1944 auch noch der Tod seiner Ehefrau dazu. All diese Vorkommnisse beeinflussten ihn offenbar unter anderem in seiner Entscheidung, sich definitiv in Chile niederzulassen[44], wo er gute Freunde gewonnen hatte und ein hohes Ansehen als Mathematiker und exzellenter Professor genoss.

### 4.2 „Die wöchentlichen 500 Stunden". Moderne Physik und Computerberechnungen/ Anfänge der Informatik

Hier bin ich nun
Hinter diesem unbequemen Pult
Verdummt vom Krach
Von fünfhundert Wochenstunden.
  Nicanor Parra, Autorretrato (1954).

Bis ins Jahr 1945 konzentrierte Professor Grandjot seine Lehrtätigkeit hauptsächlich auf das Gebiet der reinen und angewandten Mathematik. Von da an aber erweiterte er sie auch um die Physik. In den ersten Jahren des Pädagogischen Institutes, als die Lehrerausbildung drei Jahre dauerte, beschränkte sich die Unterweisung in Physik auf die experimentelle Physik. Ab 1908 wurde die Lehrerausbildung auf vier Jahre verlängert und bei dieser Gelegenheit entwickelte Poenisch den Kurs der rationalen Mechanik mit dem Ziel, den angehenden Physiklehrern eine vollständigere Ausbildung zu ermöglichen.[45] Trotz aller Sorgfalt und Strenge, mit der die Kurse der Mathematik und Physik unterrichtet wurden, blieben die Inhalte aber auf der Ebene der klassischen Themen. Daher entwickelte Grandjot 1946 den Kurs der theoretischen Physik mit der Absicht, den Schülern die moderne Physik nahezubringen. In diesem Kurs waren die wichtigsten Themen: Thermodynamik, Wellentheorie, Quantitative Mechanik und Relativitätstheorie. Als Vorlesungstext empfahl er die *Introducción a la Física Teórica* von J. Slater und N. Frank, in der der größte Teil dieser Themen behandelt wurde. Diesen Kurs hielt er bis 1962. In seinen Vorlesungen entfaltete er mit Eleganz seine vielschichtige Mathematikkultur und seine klaren und gut fundierten

---

[42] Nibaldo Bahamonde, *Discurso de incorporación como Profesor Emérito a la UMCE,* 1998. Zusätzliche Einzelheiten aus einem Gespräch mit den Autoren.

[43] J. Dieudonné, *La matemática del siglo XX,* in *La Ciencia Contemporánea,* publiziert von R. Taton, Band 4, S. 145, Edition Destino, Barcelona, 1975.

[44] Seine Einbürgerung in Chile erfolgte 1954. Seinen Wohnsitz hatte er durchgehend in Santiago de Chile.

[45] Den ersten Kurs in Rationaler Mechanik unterrichtete Gorbea im Jahre 1850 an der Fakultät der physikalischen und mathematischen Wissenschaften. Manuel Salustio Fernández, *Don Andrés Antonio Gorbea,* Anales, Mai 1861, S. 673.

Kenntnisse der modernen Physik, die er seiner-
zeit in Göttingen direkt von ihren Entdeckern
gelernt hatte: von M. Born, P. Debye, W. Hei-
senberg und anderen. Durch seine einfache
Sprache, seine Klarheit des Ausdrucks und
seine sympathische Persönlichkeit erhielten
seine Lektionen eine große Anmut und machten
in bestimmten Passagen den Eindruck, dass
gerade eben die Erschaffung der Quantitativen
Mechanik oder die Entwicklung der Relativitäts-
theorie abläuft.

1949 wurde Grandjot als Forscher an das
physikalische Institut der Universität von Chile
berufen, das erst im Jahr zuvor unter dem
Schutz des Rektorats gegründet worden war.
Zu der damaligen Zeit befand sich die physika-
lische Forschung in Chile noch in den Kinder-
schuhen. Aber es brodelte bereits unter der
Oberfläche, um dieser Disziplin den Anstoß zu
geben, damit sie Anschluss an die Fortschritte
der Kenntnisse im Rest der Welt finden könne.
In diesem Sinne schickte der Dekan der Fakul-
tät der Philosophie und Erziehung, Professor
Juan Gómez Millas, zwei junge Absolventen
des Pädagogischen Institutes und Ex-Schüler
von Grandjot, an ausgewählte europäische
Zentren, um ihre Kenntnisse der modernen
Physik zu perfektionieren: Einer zum Thema
Kosmische Strahlung, der andere zum Thema
Kristallographie. Bei ihrer Rückkehr bildeten
beide Forschungsgruppen in ihren jeweiligen
Fachthemen. In diesen Laboren entstanden
1953 die ersten internationalen Publikationen
chilenischer Physiker.[46]

Zu dieser Zeit hatte das Pädagogische Institut
bereits seinen alten Stammsitz (an der Kreu-
zung der Straßen „La Alameda" mit „Ricardo
Cumming") verlassen, in dem es mehr als ein
halbes Jahrhundert angesiedelt war und hatte
nun einen Platz in einem modernen Campus in
der Straße Macul. Dort machten die jungen
Forscher, die in ihren weißen Kitteln durch die
Gärten des Institutes promenierten, so viel Ein-
druck, dass „keiner mehr Lehrer, sondern alle
Forscher werden wollten", wie sich viel später
Raquel Martinolli, Chefin des physikalischen
Laboratoriums des Pädagogischen Institutes,
Schülerin und Nachfolgerin von Ziegler, sehn-
süchtig erinnert. 1950 gab Grandjot einen Kurs
der experimentellen Physik im Pädagogischen
Institut, in dem er besonders die Elektrizität
behandelte. Die Studenten, die als Einführung
zu diesem Thema geriebene Kämme, die an
Papierschnipsel gehalten werden, gewohnt
waren, waren sehr erstaunt, als ihr Professor
den Kurs sofort mit dem Strom direkt aus der
Steckdose des Hörsaals begann. Er verwende-
te das Lehrbuch von Pohl. Aus der dynami-
schen Elektrizität leitete er alle Konzepte und
die verwendeten Fachbegriffe ab. Die Elektro-

statik hingegen blieb als historisches Anhängsel
übrig.

Aber Grandjot beließ es nicht nur dabei, die
moderne Physik in die Universität zu bringen.
Er organisierte eine ganze Reihe von Konferen-
zen über die Relativitätstheorie und nutzte da-
bei die große Verehrung Albert Einsteins aus,
der 1955 gestorben war. An diesen nahmen er
als Hauptredner und außerdem sein Schüler
Hernán Cortés Pinto und der Ingenieur Arturo
Aldunate Phillips teil. Diese Konferenzen wur-
den sowohl in Santiago als auch in den Provin-
zen sehr gut angenommen. Das ist besonders
bemerkenswert, auch wenn die Relativitätstheo-
rie heutzutage als klassische Wissenschaft gilt.
Vielleicht ist das aber auch nur der kurzfristi-
gen Stimmung dieser Konferenzen geschuldet,
oder dem, was Grandjot – ohne viele Zahlen-
spiele – bei der feierlichen Eröffnung darbot.[47]
Er begann folgendermaßen: „Wenn ich es nun
wage, die Einstein'sche Relativitätstheorie in
einer halben Stunde und vor einem nicht-
spezialisierten Publikum zu umreißen und wenn
ich die Hoffnung hege, dabei das ganze Wesen
dieser Theorie darzulegen, einschließlich den
wichtigsten geschichtlichen Ursprüngen, den
wissenschaftlichen Grundlagen, seiner abgelei-
teten Struktur, den Schlussfolgerungen, zu de-
nen sie führt, einigen bereits eingetroffenen
Voraussagen und solchen, die noch in einer
zukünftigen Untersuchung verifiziert werden
müssen, dann fühle ich mich von der schweren
Last einer übermenschlichen Aufgabe nieder
gedrückt." Um nun aber diese Aufgabe zu lö-
sen, durchwanderte er, ähnlich wie beim Auf-
bau einer Stadt, die Bestandteile der Physik,
nämlich die Gebäude der Wärmelehre, der Me-
chanik, der Akustik und des Elektromagnetis-
mus. Er zog Vergleiche zwischen deren Fehlern
und Reparaturen, zwischen ihren Grundmauern
und Fundamenten. Er analysierte das Konzept
der „Gleichzeitigkeit" und wie Einstein den Wi-
derspruch zwischen der Physik von Newton und
den Experimenten von Michelson löste. Er kün-
digte an, dass die Details der Theorie in weite-
ren Konferenzen weiterentwickelt werden müs-
sen: „Ich selbst muss sie in späteren Tagen
darlegen […]. An einem solchen Tag werde ich
von der Äquivalenz der schweren Masse, die
wir auf der Waage messen können, und der
trägen Masse, im 2. Gesetz von Newton
verwendet wird, reden müssen. Ich werde de-
ren Gleichzeitigkeit aufzeigen müssen […]. Ich
werde die Krümmung des Raumes in 3 oder 4
Dimensionen, des Raumes von Riemann, erklä-
ren müssen; und die Identifikation dieser
Krümmung mit der universellen Gravitation; ich
werde von der Materie sprechen müssen, die sich in Interaktion mit ihrem eigenen Gravitationsfeld befindet, und davon, wie die letzten Spuren des absoluten Raumes abgeschafft werden, von privilegierten Systemen und von der innigen Vereinigung der Mechanik und des Elektromagnetismus." Das alles ist das neue und wunderschöne Gebäude, das Einstein aufgebaut hat, ein Genius, „beherzt und ehrlich, der sich nicht mit halbherzigen oder eventuell möglichen Lösungen zufrieden gibt, der stattdessen immer an der Wurzel anpackt, um die Probleme zu beseitigen, auch wenn sich dadurch in der täglichen Routine schmerzende Einschnitte ergeben." Grandjot beendete seine Konferenz mit einer Evaluation der Relativitätstheorie, indem er folgende Fragen aufwarf: Wird sie wieder außer Mode kommen? Wird sie ihre Gültigkeit verlieren? Was würde dazu Einstein sagen? Er kam dabei zu dem Schluss, dass kein Zweifel daran bestünde, dass Einstein mit Freuden für seine Theorie folgende Antwort akzeptieren würde, so wie sie auch er, Grandjot, schon formuliert hatte: „Einer Theorie kann eigentlich nichts besseres widerfahren, als dass sie von einer neuen, ausführlicheren und allgemeineren Theorie geschluckt wird."

Zusätzlich zu dem Kurs der Theoretischen Physik im Pädagogischen Institut und dem Kurs Complementos de Matemáticas Superiores, den er gleichzeitig in den Ingenieurschulen der Universität von Chile und der katholischen Universität abhielt, übernahm er 1948 einen weiteren Kurs an der Architekturschule der letztgenannten Universität. In der damaligen Zeit hatte gerade die moderne Mathematik begonnen, Einfluss auf die Architekturschulen zu nehmen. Dabei übten besonders die Mengenalgebra, die Topologie, die Gruppen- und Graphentheorie, die Themen der Analyse architektonischer Räume und urbane Symmetrien und Strukturen eine besondere Anziehung aus. Das war völlig neu für die chilenischen Architekturschulen. Grandjot bewegte sich spielerisch leicht in diesem neuen Umfeld und verband die formale Mathematik mit architektonischer Kreativität. Er wechselte von der Theorie der Verbände[48] in der Architektur zur Quantenmechanik im Pädagogischen Institut und von dort zu den Fourier-Reihen in seinen Ingenieurkursen.

Ein anderer Bereich von Aktivitäten war 1947 die Übernahme des Postens als Präsident der mathematischen Sektion der gerade gegründeten wissenschaftlichen und technologischen Forschungsdirektion der katholischen Universität (DICTUC). Seine Leistungsfähigkeit in der Arbeit war unerschöpflich. Die Vielfalt seiner Interessen ebenso. Auch für die neue Wissenschaft, die sich Stück für Stück von der Mathematik ablöste, die Anfänge der Informatik, interessierte er sich sehr. Am Ende der 50er Jahre erfasste ihn die Leidenschaft für diese neue Wissenschaft und er konstruierte einen analogen Computer mit ultrafeinen Potentiometern, um Gleichungssysteme zu lösen. „Der analoge Computer wurde im DICTUC konstruiert, um Systeme von 20 Gleichungen mit 20 Unbekannten mit einer Genauigkeit von 2% zu lösen", schrieb Grandjot in seiner Arbeit *Resolución numérica de ecuaciones algebraicas*. Ganz augenscheinlich war dieser Computer der erste, der in Chile konstruiert wurde. In diesen Zeiten standen die analoge und die digitale Informatik in Chile gleichberechtigt nebeneinander. Aber mit der Ankunft des ersten digitalen Computers, einem Lorenz ER-56, 1962 in der Universität von Chile, hatte die analoge Informatik eindeutig das Rennen verloren. Später übernahm Grandjot das Informatiklabor der katholischen Universität und gehörte damit zu den Pionieren der Informatik in Chile.

## 4.3 Ausbildung von Forschern: Der letzte Mosaikstein zu einem großartigen Lebenswerk

Nach dieser sehr intensiven akademischen Tätigkeit bis 1956, reiste Grandjot vom Januar bis Mai 1957 durch Europa und besuchte dabei Forschungszentren und –Laboratorien, ganz besonders in Italien und Deutschland. Bei seiner Rückkehr erwartete ihn bereits das durch Rektor Juan Gómez Millas gegründete Mathematische Forschungszentrum der Universität von Chile mit einem Anstellungsvertrag. Nach Aussagen von Rolando Chuaqui, dem begabtesten Schüler Grandjots, wurde mit diesem Zentrum zum ersten Mal offiziell die mathematische Forschung in Chile anerkannt.[49] Um Grandjot versammelten sich dort mehrere junge chilenische Talente, wie Legrady und andere Mathematiker, die sich voll und ganz dem Studium (dieser Disziplin) verschrieben hatten ohne nach dem Erlangen von Titeln oder Graden zu trachten. Auch ohne besondere Phantasie kann man sich gut vorstellen, wie die Beschäftigung mit Logik, Mengenlehre, Zahlentheorie, nichteuklidischer Geometrie und abstrakter Algebra sowohl Professoren als auch Studenten erfreue.

Wir haben nun ausführlich dargestellt, was für lange Umwege es bedurfte, um an diesem her-ausragenden Punkt anzulangen: von Gorbea, bei dem die Mathematik aus nützlichen Lehrsätzen für Landvermesser bestand, bis Poenisch und Tafelmacher, bei denen die Dramaturgie der Mathematik sowohl die Anwendungen für die Ingenieurswissenschaften beinhaltete als auch die Ausbildung der Lehrer für deren Tätigkeiten in den Schulen. Aber der chilenischen mathematischen Kultur fehlten noch immer die Anknüpfungspunkte, die es ihr erlaubten, aus sich heraus zu wachsen, nämlich die Ausbildung von Forschern nach ihren Prinzipien. Diese Aufgabe übernahm nun als erstes das Mathematische Forschungszentrum und später die Fakultät der Wissenschaften, beides Institute, in denen Grandjot einer der Hauptakteure war. Neben der gemeinsam abzuleistenden und offiziellen Aufgabe war es Grandjot sehr wichtig, eigene Schüler auszubilden. Im Folgenden möchten wir zwei davon besonders herausstellen, beide von großem Ansehen unter ihresgleichen, der eine in abstrakter Algebra und der andere in Logik und den Grundlagen der Mathematik. Wie wir ausgeführt haben, sind beides Disziplinen, die Grandjot gleich nach seiner Ankunft in Chile einführte und die heute durch prestigeträchtige Forschungsrichtungen im Land repräsentiert werden.

---

[48] Grandjot war ein Kenner auf diesem Gebiet und hatte das Buch von *Birkhoff*, Lattice Theory, 3. Auflage, 1967, überarbeitet. Dort hatte ihm und anderen Mathematikern Birkhoff gedankt: „Folgenden Personen bin ich zutiefst dankbar: Kirby Baker, Orrin Frink, George Grätzer, C. Grandjot, Alfred Hales, Paul Halmos, Samuel H. Holland, M. F. Janowitz, Roger Lyndon, Donald MacLaren, Richard S. Pierce, George Raney, Gian-Carlo Rota, Walter Taylor und Alan G. Waterman."

[49] R. Chuaqui, Una visión de la Comunidad Científica Nacional, CPU, Santiago, 1982, S. 12

Die abstrakte Algebra kam, wie beschrieben, sehr verspätet in die chilenischen Hörsäle. Der erste offizielle Kurs dieser „neuen Wissenschaft" wurde von Professor César Abuauad im Pädagogischen Institut der Universität von Chile gehalten, der den Lehrstuhl der Höheren Algebra von Professor Carlos Videla übernommen hatte. Abuauad war Schüler von Videla und von Grandjot. Er beendete seine Ausbildung im Pädagogischen Institut in den 30er Jahren und vertiefte sich, angeleitet von Grandjot, in die „Moderne Mathematik". Später festigte er seine Kenntnisse durch Studien in den Vereinigten Staaten im Jahre 1950. Auch wenn er nicht exakt als Forscher bezeichnet werden kann, so erlaubte ihm sein Interesse für diese Disziplin, deren Entwicklung während des 20. Jahrhunderts ganz zu verfolgen. In der Fakultät der Wissenschaften, deren Professor er viele Jahre war, erinnert man sich seiner als wichtigen Pfeiler in der Lehre der Algebra.

Der zweite von Grandjots Schülern, den wir hier erwähnen möchten, war Dr. Rolando Chuaqui, von dem auch weiter vorn im Text schon die Rede war. Chuaqui studierte eigentlich Medizin, aber während Privatstunden bei Grandjot noch während seines Studiums entdeckte er seine Liebe zur Mathematik. Chuaquis Cousin erzählte uns, dass „diese Privatstunden einmal wöchentlich, immer samstags nachmittags oder sonntags morgens stattfanden, zum Teil auch während der Ferien. Für die Erklärungen benutzte Professor Grandjot statt der Wandtafel ein Notizheft. Es wurden dabei mehrere Hefte gefüllt, von denen die meisten auch noch vorhanden sind. Die Aufgaben bestanden meist darin, ein bereits gelöstes Theorem zu beweisen, was Rolando mit großer Leichtigkeit erfüllte."[50] Sein Cousin informierte uns auch, dass auf Grund der formalistischen Einstellung Grandjots die Logik das erste und längste abgehandelte Thema war. Danach wandten sie sich der nichteuklidischen Geometrie zu, die Grandjot in axiomatischer Weise erklärte und zum Schluss kam etwas Zahlentheorie und Moderne Algebra. Nach dem Abschluss als Chirurg ging Rolando Chuaqui nach Berkeley, wo er seinen Doktortitel in Mathematik mit einer Arbeit auf dem Gebiet der Grundlagen der Wahrscheinlichkeitsrechnung erhielt.[51] Nach seiner Rückkehr nach Chile war er sehr produktiv auf dem Feld der Logik als Forscher tätig, was sowohl national als auch international anerkannt wurde.

Grandjot blieb bis zum Ende des Jahres 1959 im Zentrum der mathematischen Forschung. Von 1960 bis 1963 bekleidete er im Institut der Physik und Mathematik der Ingenieursschule der Universität von Chile das Amt eines Professors und Forschers und hielt Seminare zu verschiedenen Feldern dieser Disziplinen ab. 1965 und 1966 erfreute sich die Fakultät der Wissenschaften seiner Dienste. Dieses Institut war Nachfolger des Mathematischen Forschungszentrums geworden mit dem Plan, unter anderem Bachelor- und Magistertitel in den Wissenschaften zu vergeben.[52] Dem Vorbild folgten bald weitere Institutionen in allen chilenischen Universitäten. Gipfelpunkt dieses Prozesses war 1967 die Einsetzung der Nationalen Kommission für wissenschaftliche Forschung und Technologie, die direkt der Regierung unterstellt war. Bei dieser Entwicklung war Professor Juan Gómez Millas, ehemaliger Rektor der Universität von Chile und ehemaliger Bildungsminister die treibende Kraft. „Von ihm kamen die fortschrittlichsten Impulse, die das wissenschaftliche Chile je erlebt hat."[53] Auch Grandjot spielte hierbei eine wichtige Rolle als Professor und Berater, vor allem bei der Ausbildung und bei der Auswahl von jungen Menschen, die ausgezeichnet wurden, indem man sie zu den anspruchsvollsten wissenschaftlichen Zentren der ganzen Welt schickte. Wir denken, dass sich Grandjot dabei in Verantwortung sah, die noch unerledigten Aufgaben abzuarbeiten, die ihm seine Vorgänger hinterlassen hatten, besonders die, die Mathematik in Chile zu verbreiten, die nicht nur eine „nützliche" Wissenschaft oder eine Ansammlung von Wissen, sondern auch einen „kreativen Prozess" in sich darstellt.

---

[50] B. Chuaqui, *Rolando Chuaqui en la Escuela de Medicina*, 2001, unveröffentlicht.
[51] Es ist bekannt, dass Grandjot ein Empfehlungsschreiben für Rolando schrieb, dass offensichtlich großen Einfluss hatte, denn Rolando war in Chile nicht als Student der Mathematik eingeschrieben. Es wäre sicherlich interessant, in den Wortlaut dieses Empfehlungsschreibens Einsicht zu bekommen.
[52] Die Fakultät der Wissenschaften der Universität von Chile wurde 1965 ins Leben gerufen. Dies geschah aber erst nach langen und intensiven akademischen Diskussionen, was zeigt, mit welchen Schwierigkeiten die Institutionalisierung der Wissenschaften in Chile zu kämpfen hatte. Gleichzeitig (1965) bekundete der höheren Rates der Katholischen Universität, laut Prof. Salinas, „dass eine solche Fakultät nur „zu einer Fabrik zur Produktion von Wissenschaftlern" werden würde und dass so etwas das Land verschlingen würde". Sie meinten auch, dass diese nur ein Rettungsring für mittelmäßige Schüler wäre und, zu guter Letzt, eine solche Fakultät zu viel Geld verschlingen würde." (Zitiert durch Raúl Sáez, in *Universidad, Ciencia y Desarrollo (Universität, Wissenschaft und Entwicklung)*, Hombres del siglo XX, Band II, Dolmén Ediciones, 1994, S. 1213.)

## 4.4 Seine letzten Jahre: Krankheit und Vergessen

Am 17.11.1966 erließ das Rektorat der Universität von Chile folgendes Dekret 9397:

„Vorlage: das Beschlossene im Abschnitt 3°, Absatz 2° des DFL 338 von 1960 und die *unerlässliche Notwendigkeit*, dass Don Carlos Grandjot Reins am Mathematikerkongress in Moskau teilnimmt.

Beschluss: Don Carlos Grandjot Reins (Rol C.G. 37262), Ordentlicher Professor der Mathematik in der Zentralen Abteilung der mathematischen- und Naturwissenschaften, wird für einen Monat mit dem Zweck abgeordnet, am Mathematikerkongress in Moskau, Russland, teilzunehmen. Die Kommission ist sehr erfreut darüber." (Hervorhebung in kursiv durch die Autoren).

Dieses Dekret ist von Rektor Eugenio González unterschrieben worden. Die Aussage „unerlässliche Notwendigkeit" zeugt ohne Zweifel vom Interesse der Universitätsbehörden und der Regierung, die neuen Institutionen und deren wissenschaftlichen Aktivitäten mit den entsprechenden Institutionen in Amerika und Europa zu verknüpfen. Entsprechend seinem Lebenslauf war Grandjot dafür besonders geeignet. Die Kommission verlängerte die Frist bis zum 28. Februar 1967.

Nach Beendigung seiner Abordnung kam Grandjot Ende Februar 1967 wieder nach Chile zurück. Während ihrer gemeinsamen Reise nach Moskau und Deutschland war im November 1966 seine zweite Frau, Margarete Rixmann Hövener, gestorben.

Am 1. März 1967, nach 38 Jahren ununterbrochenem Dienst an der chilenischen Wissenschaft und Ausbildung ging Carlos Grandjot in den Ruhestand. Als unermüdlicher Arbeiter machte er jedoch auch danach noch als außerordentlicher Professor an der Universität von Chile weiter. Leider erlitt er am Ende desselben Jahres einen Schlaganfall, der seinen akademischen Aktivitäten ein jähes und endgültiges Ende setzte.

Es folgten Jahre der Einsamkeit und des Vergessens. Eine äußerst traurige Periode seines Lebens. Die Aufmerksamkeit seiner Kollegen, die ihn hätten besuchen können, wurde abgelenkt durch die universitären Reformen, durch die politischen Kämpfe, die in dieser Zeit ausgetragen wurden, durch den Militärputsch im September 1973, der die Universitäten beeinträchtigte. Im Oktober 1973 holte ihn seine Tochter Sigrid nach Concepción, wo sie als Mathematiklehrerin arbeitete. Dort verbrachte er die letzten Jahre seines Lebens im Deutschen Altersheim, in dem ihn seine Tochter jeden Nachmittag besuchte, um mit ihm seine undeutlichen Erinnerungen aufzufrischen: an die wissenschaftliche Gesellschaft, an seine Universitätskollegen, seine Schüler und alle Institutionen, denen er seine besten Jahre gewidmet hatte.[54]

Er starb am 5. Oktober 1979 und wurde im Hauptfriedhof von Concepción begraben.

## Nachwort

Durch einen Telefonanruf von Señora Marta Gutmann, einer Freundin Sigrids, erhielten die Autoren die Nachricht, dass Sigrid Grandjot Fritsche am Montag, den 2. Dezember 2002 in Rancagua (Chile) gestorben war, zu einem Zeitpunkt, als dieser Bericht kurz vor seiner Fertigstellung war. Sie wurde gemäß ihrem Wunsch an der Seite ihres Vaters beerdigt. Mit ihr starb der einzige Nachfahre der Familie Grandjot in Chile.

---

[53] Héctor Croxatto, *Una visión de la Comunidad Científica Nacional (Eine Zukunftsvision einer nationalen wissenschaftlichen Gemeinschaft)*, CPU, 1982.

[54] Aus Gesprächen zwischen Sigrid und den Autoren.

**5 Zeittafel**

| | |
|---|---|
| 1900 | Geburt am 23. August in Frankenberg, Deutschland. |
| 1919 | Aufnahme in der Universität von Göttingen. |
| 1922 | Erhalt des Doktortitels in Philosophie mit Erwähnung der Mathematik, Universität von Göttingen. |
| 1925 | Er erhält den Titel Privatdozent und beginnt seine Vorlesungstätigkeit an der Universität von Göttingen. |
| 1926 | Heirat mit Gertrudis Fritsche. |
| 1928 | Stipendium der Rockefeller-Stiftung in Paris. |
| 1929 | Geburt seiner Tochter Sigrid in Paris, 13. Februar. |
| | Am 9. April Unterschrift des Vertrages mit der Chilenischen Regierung über 2 Jahre, wenige Tage danach Einschiffung nach Chile. |
| | Ankunft in Chile am 1. Mai. |
| | Beginn seiner Lehrtätigkeit im Pädagogischen Institut. |
| | Im September erreicht seine Frau mit der kleinen Tochter Chile. |
| 1930 | 30. Dezember: Gründung des Chilenischen Institutes der Wissenschaften. Grandjot ist einer seiner Gründerväter. |
| 1931 | Erste Reise nach Deutschland. |
| 1933 | Beginn der Vorlesungen in der Fakultät für Ingenieure der Pontificia Universidad Católica (=PUC, katholische Universität). |
| 1938 | Zweite Reise nach Deutschland. |
| 1940 | Grandjot schreibt seine „Algebra Abstracta". |
| 1941 | Reise durch Bolivien und Hochperu. |
| 1944 | Tod seiner Ehefrau Gertrudis Fritsche. |
| 1945 | Professur in der Ingenieursschule der Universität von Chile. |
| | Eheschließung mit seiner zweiten Frau Margarete Rixmann Hövener |
| 1947 | Ernennung zum Chef der neu gegründeten Mathematischen Abteilung der DICTUC der PUC. |
| 1952 | Ratsmitglied der chilenisch-deutschen Liga (bis 1959). |
| 1953 | Professor der Architekturschule der PUC. |
| | Präsident und Gründungsvater der Mathematischen Chilenischen Gesellschaft. |
| 1954 | Erhalt der chilenischen Staatsbürgerschaft. |
| 1962 | Ernennung zum Chef des Laboratoriums für elektronische Rechner der PUC. |
| 1966 | Teilnahme am Mathematischen Kongress in Moskau. |
| | Tod seiner zweiten Ehefrau in Deutschland. |
| 1967 | Am 1. März Eintritt in den Ruhestand nach 38 Jahren Dienst in Chile. |
| | Am Jahresende erleidet er einen Schlaganfall, von dem er sich nicht mehr erholt. |
| 1973 | Im Oktober holt ihn seine Tochter in das Deutsche Altersheim in Concepción. |
| 1979 | Er stirbt am 5. Oktober in Concepción. |

**Schriftenverzeichnis von Carlos Grandjot**

**Literatur**

**[A] Bücher und Artikel**

**[B] Zeitschriften und Zeitungen**

**[C] Interviews und Berichte:**

César Abuauad, Sigrid Grandjot, Benedicto Chuaqui, Guacolda Antoine, Nibaldo Bahamonde.

Claudio Gutiérrez,          Flavio Gutiérrez,
Universidad de Chile,       Universidad de Valparaíso,
Chile                       Chile

Anmerkungen zu dieser deutschsprachigen Version des Artikels können an Christoph Lamm, christoph.lamm@web.de, geschickt werden.